\documentclass[12pt]{article}
\usepackage{amsmath, amsthm, amssymb}
\numberwithin{equation}{section}
\begin{document}

\author{Ajai Choudhry and Arman Shamsi Zargar}
\title{Pairs of equiperimeter and equiareal triangles\\
whose sides are perfect squares} 

\date{}
\maketitle

\begin{abstract}
In this paper we consider the problem of finding pairs of triangles whose sides are perfect squares of integers, and which have a common perimeter and common area. We find two such pairs of triangles, and  prove that there exist infinitely many pairs of triangles with the specified properties.
\end{abstract}

\maketitle

Mathematics Subject Classification 2020: 11D41

Keywords: equiperimeter triangles; equiareal triangles; triangles with squared sides. 

\section{Introduction}
Ever since the discovery of right-angled triangles with integer sides, there has been considerable interest in finding triangles as well as polygons with certain geometric properties and all of whose sides are given by integers. Several mathematicians have also considered diophantine problems pertaining to a pair of triangles or other geometric objects (see for instance, \cite{AT}, \cite{BG}, \cite{DJM}, \cite{LM}, \cite{LZ}, \cite{SU}, \cite{Zh}).   
 Considerable attention has been given to the problem of finding two triangles with  the same  perimeter and the same area and such that all the  sides and the common area of the two triangles are given by integers (see \cite{Aa}, \cite{Br}, \cite{Lu}, \cite{Yi}). In fact, Choudhry \cite{Ch} has described a method of generating an arbitrarily large number of scalene rational triangles with the same perimeter and the same area. 

In this paper we consider the problem of finding two triangles with the same perimeter and the same area and such that all the sides of the two triangles are perfect squares of integers.

\section{Equiperimeter and equiareal triangles with squared sides}
Let us consider two triangles whose sides are $a^2, b^2, c^2$ and $d^2, e^2, f^2$ where $a, b, \ldots, f$, are all integers. The two triangles will have the same perimeter if the integers $a, b, \ldots, f$, satisfy the following condition:
\begin{equation}
a^2+b^2+c^2=d^2+e^2+f^2. \label{condperim}
\end{equation}
Further, on using  Heron's well-known formula for the area of a triangle, the condition that the two triangles will have the same area may be written as follows:
\begin{multline}
(a^2+b^2+c^2)(a^2+b^2-c^2)(a^2-b^2+c^2)(-a^2+b^2+c^2)\\
=(d^2+e^2+f^2)(d^2+e^2-f^2)(d^2-e^2+f^2)(-d^2+e^2+f^2). \label{condarea}
\end{multline}

Since both equations \eqref{condperim} and \eqref{condarea} are homogeneous, any rational solution of these equations will yield, on appropriate scaling, a solution in integers. It therefore suffices to obtain rational solutions of \eqref{condperim} and \eqref{condarea}. 

We note that, in view of the condition \eqref{condperim}, equation \eqref{condarea} reduces to 
\begin{multline}
(a^2+b^2-c^2)(a^2-b^2+c^2)(-a^2+b^2+c^2)\\
=(d^2+e^2-f^2)(d^2-e^2+f^2)(-d^2+e^2+f^2). \label{condarea2}
\end{multline}

To solve the simultaneous diophantine equations \eqref{condperim} and \eqref{condarea2}, we write
\begin{equation}
\begin{aligned}
a &= pu + q + r, \quad & b&= qu - p - r, \quad & c &= ru - p + q,\\
d &= pu - q - r, \quad & e &= qu + p + r, \quad & f &= ru + p - q,
\end{aligned}
\label{subs1}
\end{equation}
where $p, q, r$ and $u$ are arbitrary parameters.
 
With these values of $a, b, \ldots, f$, it is readily verified that equation \eqref{condperim} is identically satisfied. Further, equation \eqref{condarea2} reduces to
\begin{multline}
16u(u - 1)(u + 1)(q + r)(p + r)(p - q)\\
\times \{(p^3 + p^2q - p^2r+ pq^2 - 2pqr + pr^2 + q^3 - q^2r + qr^2 - r^3)u^2\\
+2p^2q - 2p^2r + 2pq^2 + 12pqr + 2pr^2 - 2q^2r + 2qr^2\}=0. \label{condarea3}
\end{multline}

We note that when the parameters $p, q, r$ and $u$ satisfy the condition, 
\begin{equation}
u(u - 1)(u + 1)(q + r)(p + r)(p - q)(p+q-r) = 0, \label{trivcond}
\end{equation}
we get  a trivial solution of equations \eqref{condperim} and \eqref{condarea}. Accordingly,  we must find rational numbers $p, q, r$ and $u$ such that the above condition is not satisfied 
and the last factor of \eqref{condarea3} becomes 0. Thus  there must exist an integer  $t$ such that 
\begin{multline}
t^2=-(p^3 + p^2q - p^2r + pq^2 - 2pqr + pr^2 + q^3 - q^2r + qr^2 - r^3)\\
\times (2p^2q - 2p^2r + 2pq^2 + 12pqr + 2pr^2 - 2q^2r + 2qr^2). \label{condt}
\end{multline}
 
Further,  to obtain actual triangles with a common perimeter and common area, we need to find solutions of the simultaneous equations \eqref{condperim} and \eqref{condarea} such that  the sides $a^2, b^2, c^2$ and $d^2, e^2, f^2$ of the two triangles satisfy the triangle inequalities. A computer search in the range $|p| + |q| +|r|  \leq 700$ yielded two sets of values of $p, q, r$ such that all conditions are satisfied.

The first pair of triangles is obtained by taking $(p, q, r)=(14, -27, -25)$ when we get two triangles with sides 
\begin{equation}
 661^2, 1498^2, 1515^2 \quad {\rm and} \quad  921^2, 1310^2, 1553^2, \label{firstpair}
\end{equation}
 such that the two triangles have a common perimeter and common area.

A second example is obtained with $(p, q, r)=(46, 73, 371)$ when again we get two triangles with a common perimeter and common area, the sides of the two triangles being 
\begin{equation}
71297^2, 77895^2, 97154^2 \quad {\rm and} \quad  67005^2, 81926^2, 96893^2. \label{secpair}
\end{equation}

We will now show that there exist infinitely many pairs of triangles with a common perimeter and common area and such that all the sides of both the triangles are given by perfect squares of integers. 

On writing 
\begin{equation}
q=-p(1-mx),\quad r=px,\quad t=p^3xy, \label{valqrt}
\end{equation}
the condition \eqref{condt} may be written as follows:
\begin{equation}
\begin{aligned}
y^2& = 2(m - 1)^2(m^2 + 1)mx^4 - 2m(m - 1)(m^3 + 10m^2 + m + 8)x^3\\
& \quad \quad  + (6m^4 + 50m^3 - 18m^2 + 14m - 16)x^2 - 8m^2(m + 8)x\\
 & \quad \quad+ 4m(m + 8). 
\end{aligned}
\label{quartec}
\end{equation}

Any rational solution of equation \eqref{quartec} yields, on using the relations \eqref{valqrt}, a solution of  the equation \eqref{condt}. Moreover, the solution $(p, q, r)=(14, -27, -25)$  of condition \eqref{condt} corresponds to the solution 
\[(m, x, y)=(13/25, -25/14, -339/245)\]
of equation \eqref{quartec}. 

Accordingly, we fix $m=13/25$ in equation \eqref{quartec} when we get 
\begin{equation}
\begin{aligned}
y^2&=(2972736/9765625)x^4 + (55402464/9765625)x^3\\
 & \quad \quad- (2389884/390625)x^2 - (287976/15625)x + 11076/625.
\end{aligned}
\label{quartecex1}
\end{equation}

Now equation \eqref{quartecex1} represents a quartic model of an elliptic curve, and the birational transformation given by
\begin{equation}
\begin{aligned}
x &= -25(2048250X + 4633Y- 1188723263160)\\
&\quad \quad \times(140122182X - 23657Y-532179246194760)^{-1},\\
y&=5061168(137842X^3 - 739070924100X^2 - 68921Y^2\\
 &\quad \quad - 182461764476160Y + 3026923795437314317841600)\\
&\quad \quad \times\{5(140122182X - 23657Y-532179246194760)\}^{-2}
\end{aligned}
\label{biratxy}
\end{equation}
and
\begin{equation}
\begin{aligned}
X& =4(1311254697763x^2 + 2523717076250x + 1621155375000y\\
 &\quad \quad - 3370847328125)(577x-2825)^{-2},\\
Y&=60734016(223993723304x^3 + 1368674740750x^2\\
 &\quad \quad  + 364901515625xy- 2217844262500x + 133349609375y\\
 &\quad \quad- 372517031250)(577x-2825)^{-3}\\
\end{aligned}
\label{biratXY}
\end{equation}
reduces the quartic curve \eqref{quartecex1} to the Weierstrass form given by
\begin{equation}
Y^2 = X^3 - 21151030877616X + 31685265497576201600. \label{cubicecex1}
\end{equation}

Using the software SAGE, we quickly found the rank on the elliptic curve \eqref{cubicecex1} to be 2, with  the generators being 
\begin{equation*}
\begin{aligned}
&(6008706700/1681,  91230882238080/68921), \\
&{\rm and} \quad (7840706250956/1168561, -17496345598032878080/1263214441).
\end{aligned}
\end{equation*}

Now we already know a rational point $(x, y)=(-25/14, -339/245)$ on the quartic curve \eqref{quartecex1}, and corresponding to this rational point, we readily find a  point $P$ on the cubic curve \eqref{cubicecex1}, the coordinates of $P$ being as follows:
\begin{equation}
(-7450305309428/4661281, -78862809542759294976/10063705679).
\end{equation}

It is clear from the manner in which we have obtained the elliptic curve \eqref{cubicecex1} that corresponding to each rational point on the curve \eqref{cubicecex1}, there exists a rational solution of the simultaneous equations \eqref{condperim} and \eqref{condarea}. The point $P$ on the  curve \eqref{cubicecex1} corresponds to the solution $(p,q,r)=(14, -27, -25)$  of equation \eqref{condt} and yields the two triangles mentioned at  \eqref{firstpair}.

Since the curve \eqref{cubicecex1} has positive rank, there are infinitely many rational points on the elliptic curve \eqref{cubicecex1} and it  follows from   a theorem of Poincare and Hurwitz \cite[Satz 11, p. 78]{Sk} that  there exist infinitely many rational points on the curve \eqref{cubicecex1} in the neighbourhood of the point $P$, and these points will yield  infinitely many rational solutions  of the simultaneous equations \eqref{condperim} and \eqref{condarea} such that $a^2, b^2, c^2$ and $d^2, e^2, f^2$ satisfy the triangle inequalities, and we can thus  generate infinitely many pairs of triangles whose sides are perfect squares of integers, and which have a common perimeter and common area.

\section{Concluding remarks}

The aforementioned theorem of Poincare and Hurwitz does not give any method of generating  infinitely many rational points on the curve \eqref{cubicecex1} in the neighbourhood of the point $P$, and hence is of no help in actually finding pairs of triangles whose sides are perfect squares of integers, and which have a common perimeter and common area. It would be of interest to find an effective method of generating infinitely many pairs of triangles with the desired properties.

Finally we note that while the triangles we have obtained have a common perimeter and common area, the area of our triangles is not an integer. It is very unlikely that there exist pairs of triangles  whose sides are perfect squares of integers, and which have a common perimeter and common area which is also an integer.

\noindent Postal Address: Ajai Choudhry, 13/4 A Clay Square, \\
\noindent Lucknow - 226001, India\\
\noindent E-mail: ajaic203@yahoo.com

\medskip

\noindent Postal Address: Arman Shamsi Zargar, Department of Mathematics \\
\noindent and Applications, University of Mohaghegh Ardabili,\\
\noindent  Ardabil, IRAN\\
\noindent E-mail: zargar@uma.ac.ir

\end{document}